# SMALL DEVIATIONS OF GENERAL LÉVY PROCESSES

By Frank Aurzada[1] and Steffen Dereich[2]

*Technische Universität Berlin*

We study the small deviation problem $\log \mathbb{P}(\sup_{t \in [0,1]} |X_t| \leq \varepsilon)$, as $\varepsilon \to 0$, for general Lévy processes $X$. The techniques enable us to determine the asymptotic rate for general real-valued Lévy processes, which we demonstrate with many examples.

As a particular consequence, we show that a Lévy process with nonvanishing Gaussian component has the same (strong) asymptotic small deviation rate as the corresponding Brownian motion.

## 1. Introduction and results.

1.1. *Motivation and notation.* The small deviation problem for a stochastic process $X = (X_t)_{t \in [0,1]}$—also called the small ball problem—consists in determining the probability

$$\mathbb{P}\left(\sup_{t \in [0,1]} |X_t| \leq \varepsilon\right) \qquad \text{as } \varepsilon \to 0.$$

One can also consider other norms, but, in this article, we concentrate on the supremum norm, which is denoted by $\|\cdot\|$. There has been a lot of interest in small deviation problems in recent years, which is due to the many connections to other questions, such as the law of the iterated logarithm of Chung type, strong limit laws in statistics, metric entropy properties of linear operators, quantization and several other approximation quantities for stochastic processes (see the surveys [10] and [9] and the bibliography [11]).

Typically, one cannot determine the above probability, even asymptotically, except for a very few examples (such as, e.g., Brownian motion).

Received May 2008; revised November 2008.

[1]Supported by the DFG Research Center MATHEON "Mathematics for key technologies" in Berlin.

[2]Partly supported by a DFG Research Fellowship.

Therefore, one concentrates on the asymptotic rate of the logarithm of that quantity,

$$-\log \mathbb{P}\Big(\sup_{t\in[0,1]}|X_t|\le\varepsilon\Big) \qquad \text{as } \varepsilon\to 0. \tag{1}$$

Even this simplified problem is a difficult issue if one aims to solve it for *a whole class of processes.* Thus far, this has only been possible for a large subclass of Gaussian processes (see the approach in [7], completed in [8]). With the framework presented in this article however, we are able to determine the asymptotic rate of the quantity (1) for general real-valued Lévy processes.

Let $X = (X_t)_{t\in[0,1]}$ denote a Lévy process. It is characterized by independent and stationary increments, $X_0 = 0$, stochastic continuity and cadlag paths; see [3, 18]. Due to the Lévy–Khintchine formula, the characteristic function of each marginal $X_t$ $(t\in[0,1])$ admits the representation

$$\mathbb{E}e^{iuX_t} = e^{-t\psi(u)}, \tag{2}$$

where

$$\psi(u) = \frac{\sigma^2}{2}u^2 - ibu + \int_{\mathbb{R}\setminus\{0\}}(1 - e^{iux} + \mathbb{1}_{\{|x|\le 1\}}iux)\nu(dx)$$

for parameters $\sigma^2\in[0,\infty)$, $b\in\mathbb{R}$ and a positive measure $\nu$ on $\mathbb{R}\setminus\{0\}$, called *Lévy measure*, satisfying

$$\int_{\mathbb{R}\setminus\{0\}} 1\wedge x^2\,\nu(dx) < \infty. \tag{3}$$

On the other hand, for a given triplet $(\nu,\sigma^2,b)$, there exists a Lévy process $X$ such that (2) is valid, and its distribution is uniquely characterized by the latter triplet. We call the corresponding process a $(\nu,\sigma^2,b)$-*Lévy process.* In order to avoid pathological cases, we always assume that the Lévy process is nondeterministic.

We recall that, for Brownian motion $B$, $\sigma > 0$ and $b\in\mathbb{R}$,

$$-\log \mathbb{P}\Big(\sup_{t\in[0,1]}|\sigma B_t + bt|\le\varepsilon\Big) \sim \frac{\pi^2}{8}\sigma^2\varepsilon^{-2}, \tag{4}$$

where $\sim$ means strong asymptotic equivalence. This process corresponds to a $(0,\sigma^2,b)$-Lévy process; see, for example, Section 7.3 in [12] for historical remarks on this result.

Apart from this special case, the rate in (1) is already known for certain stable Lévy processes. Namely, for any *strictly* $\alpha$-stable Lévy process [i.e., a Lévy process satisfying the scaling property $\mathcal{L}(X_t) = \mathcal{L}(t^{1/\alpha}X_1)$ for all $t\ge 0$] such that $|X|$ is not a subordinator (an increasing process, see below), we have

$$-\log \mathbb{P}\Big(\sup_{t\in[0,1]} |X_t| \le \varepsilon\Big) \sim K \varepsilon^{-\alpha} \tag{5}$$

for some constant $K > 0$ (see [3], page 220, or [5, 16] and [23]). We refer to [13] and [21] for an overview of further results for symmetric stable Lévy processes under various kinds of other norms.

On the other hand, if $|X|$ is a stable subordinator [then, necessarily, $0 < \alpha < 1$ and its characteristic function is given by (8) below], then, for some explicitly known $K' > 0$, as $\varepsilon \to 0$,

$$-\log \mathbb{P}\Big(\sup_{t\in[0,1]} |X_t| \le \varepsilon\Big) = -\log \mathbb{P}(|X_1| \le \varepsilon) \sim K' \varepsilon^{-\alpha/(1-\alpha)}. \tag{6}$$

In the case where $X$ is a subordinator, it is easy to determine the small deviation rate since, here, $\sup_{t\in[0,1]} |X_t| = X_1$. It thus suffices to look at the distribution of $X_1$ in a neighborhood of the origin. If $b > \int_0^1 x\nu(dx)$, then the probability of a small deviation is zero for sufficiently small $\varepsilon$, whereas, in the case $b = \int_0^1 x\nu(dx)$, one can determine the asymptotics of (1) via de Bruijn's Tauberian theorem ([4], Theorem 4.12.9) from the asymptotics of

$$-\log \mathbb{E}e^{-\lambda X_1} = -\int_0^\infty (e^{-\lambda x} - 1)\nu(dx) \qquad \text{as } \lambda \to \infty, \tag{7}$$

provided that the latter expression is regularly varying.

Let us introduce further notation. If (2) is true for

$$\psi(u) = \frac{\sigma^2}{2}u^2 + \int_{\mathbb{R}\setminus\{0\}} (1 - e^{iux} + iux)\nu(dx)$$

for a measure $\nu$ with $\int |x| \wedge x^2 \nu(dx) < \infty$, then we call $X$ a $(\nu, \sigma^2)$-*Lévy martingale*. It is a martingale in the usual sense. Furthermore, we say that a $(\nu, \sigma^2, b)$-Lévy processes is *of type (I)* if

$$\int_{-1}^1 |x|\nu(dx) < \infty \quad \text{and} \quad \sigma^2 = 0.$$

Finally, a $(\nu, \sigma^2, b)$-Lévy process is called a *subordinator* if it is almost surely increasing. Recall (see [18]) that this is the case if and only if $\sigma^2 = 0$, $\nu$ is concentrated on the positive real line and satisfies

$$\int_0^1 x\nu(dx) < \infty \quad \text{and} \quad b \ge \int_0^1 x\nu(dx).$$

We use the following notation for, respectively, strong and weak asymptotics. We write $f \lesssim g$ if $\limsup f/g \le 1$; $f \gtrsim g$ is defined analogously.

Further, $f \sim g$ means that $\lim f/g = 1$. We also use $f \preceq g$ (or $g \succeq f$) if $\limsup f/g < \infty$ and $f \approx g$ if $0 < \liminf f/g \leq \limsup f/g < \infty$.

This paper is organized as follows. In Section 1.2, we review results of Simon, who studied the question when the problem (1) actually makes sense for Lévy processes. Section 1.3 contains our main results, which are illustrated by several examples in Section 2. The proofs are postponed to Sections 3 and 4 (proofs of the main results) and Section 5 (proof of the explicit rates in the examples).

1.2. *The small deviation property.* In [19], the following question was studied: for which Lévy processes does the small deviation problem make sense? Namely, one says that a stochastic process $X = (X_t)_{t \in [0,1]}$ possesses the small deviation property if

$$\mathbb{P}\Big(\sup_{t \in [0,1]} |X_t| \leq \varepsilon\Big) > 0 \qquad \text{for all } \varepsilon > 0.$$

Simon investigated this property for $\mathbb{R}^d$-valued Lévy processes. For real-valued Lévy processes, it reduces to the following, easily verifiable, equivalent characterization [19].

PROPOSITION 1.1. *A $(\nu, \sigma^2, b)$-Lévy process $X$ possesses the small deviation property if and only if it is not of type* (I) *or if it is of type* (I) *and, for $c := b - \int_{|x| \leq 1} x \nu(dx)$, we have:*

- *$c = 0$, or*
- *$c > 0$ and $\nu\{-\varepsilon \leq x < 0\} \neq 0$ for all $\varepsilon > 0$, or*
- *$c < 0$ and $\nu\{0 < x \leq \varepsilon\} \neq 0$ for all $\varepsilon > 0$.*

Let us visualize this fact with a simple example.

EXAMPLE 1.2. Let us consider an $\alpha$-stable subordinator with drift: $X_t + \mu t$, where $X$ has the characteristic function:

$$(8) \qquad \mathbb{E}e^{izX_t} = \exp\Big( t \int_0^\infty (e^{izx} - 1)\frac{dx}{x^{1+\alpha}} \Big)$$

and, by Proposition 1.1, $X_t + \mu t$ possesses the small deviation property if and only if $\mu \leq 0$, that is, if there is a nonpositive drift. Clearly, if there is a positive drift, then, already, the drift term makes the process leave the interval $[-\varepsilon, \varepsilon]$ almost surely for $\varepsilon < \mu$. For $\mu = 0$, relation (6) holds, whereas the rate of (1) was previously unknown for $\mu < 0$. We come back to this case in Example 2.3.

Let us comment on some related work. Ishikawa [6] generalizes the results from Simon [19] to other types of stochastic processes. A related result for scaled Poisson processes connected to the Strassen law is shown in [1]. Further, similar results are obtained in [22] for symmetric stable Lévy processes with $1 < \alpha < 2$. Moreover, Simon [20] investigated the small deviation problem for Lévy processes under $p$-variation norm in contrast to the present considerations. Finally, the works [14] and [15] solve the small deviation problem for certain specific Lévy processes, namely, Lévy processes that arise from subordination to Brownian motion. We come back to this relation in Examples 2.12 and 2.13.

1.3. *Main results.*  In this section, we show how to obtain estimates for (1) for a general real-valued Lévy process $X$ with triplet $(\nu, \sigma^2, b)$. Our approach is based on two results which we now state.

PROPOSITION 1.3.  *Let $(X_t)$ be a Lévy process whose Lévy measure has support in $\{|x| \leq 1\}$ and assume that $u^* \in \mathbb{R}$ is a solution of*

$$(9) \qquad\qquad \Lambda'(u^*) = 0,$$

*where*

$$(10) \qquad \Lambda(u) := \frac{1}{2}\sigma^2 u^2 + bu + \int [e^{ux} - 1 - ux]\nu(dx)$$

*denotes the logarithmic moment generating function of $X_1$. Then, the Esscher transform $\mathbb{Q}$ given by*

$$\frac{d\mathbb{Q}}{d\mathbb{P}} = e^{u^* X_1 - \Lambda(u^*)}$$

*is a probability measure such that, for all $\varepsilon > 0$,*

$$(11) \quad e^{\Lambda(u^*) - \varepsilon|u^*|}\mathbb{Q}(\|X\| \leq \varepsilon) \leq \mathbb{P}(\|X\| \leq \varepsilon) \leq e^{\Lambda(u^*) + \varepsilon|u^*|}\mathbb{Q}(\|X\| \leq \varepsilon)$$

*and $X$ is a $(e^{u^* x} \cdot \nu(dx), \sigma^2)$-Lévy martingale under $\mathbb{Q}$.*

We remark that Proposition 1.3 can be shown in a more general context, in particular, for a broad class of noncompactly supported Lévy measures and in the multidimensional setting. However, the current formulation is sufficient for our purposes.

This result enables us to transform a Lévy process with compactly supported Lévy measure into a Lévy martingale. It turns out that we find an appropriate Esscher transform [i.e., (9) has a unique solution] in all cases needed. To be more precise, Proposition 1.3 can be applied if $X$ (or $-X$) is not a subordinator and $X$ possesses the small deviation property. The proof of this fact (formulated as Lemma 3.2) and the proof of Proposition 1.3 are given in Section 3.

We proceed with stating the second key result, which is proved in Section 4. It treats the case of a Lévy martingale.

PROPOSITION 1.4. *Let $\varepsilon > 0$ and denote by $X$ a $(\nu, \sigma^2)$-Lévy martingale with $\nu$ supported on $[-\varepsilon, \varepsilon]$. Then,*

$$\mathbb{P}(\|X\| \leq 3\varepsilon) \geq e^{-10F(\varepsilon)-3} \quad and \quad \mathbb{P}(\|X\| \leq \varepsilon/2) \leq e^{F(\varepsilon)/12+1},$$

*where*

$$(12) \qquad F(\varepsilon) := \frac{1}{\varepsilon^2}\left[\sigma^2 + \int_{-\varepsilon}^{\varepsilon} x^2 \nu(dx)\right].$$

We now outline how Propositions 1.3 and 1.4 lead to good estimates for the logarithmic small ball probabilities (1). Let $X$ be an arbitrary $(\nu, \sigma^2, b)$-Lévy process. We fix $\varepsilon > 0$ and denote by $\bar{\mathbb{P}}^{\varepsilon}$ the conditional probability of $\mathbb{P}$ given that $X$ has no jumps that are larger than $\varepsilon$. Under $\bar{\mathbb{P}}^{\varepsilon}$, $X$ is a Lévy process with triplet $(\nu|_{[-\varepsilon,\varepsilon]}, \sigma^2, b_\varepsilon)$ with $b_\varepsilon := b - \int_{[-1,1]\setminus[-\varepsilon,\varepsilon]} x\nu(dx)$. Thus, under $\bar{\mathbb{P}}^{\varepsilon}$, $X_1$ has the logarithmic moment generating function

$$(13) \qquad \Lambda_\varepsilon(u) = \frac{1}{2}\sigma^2 u^2 + b_\varepsilon u + \int_{-\varepsilon}^{\varepsilon}[e^{ux} - 1 - ux]\nu(dx).$$

Assuming that there exists a solution $u_\varepsilon$ to the equation $\Lambda'_\varepsilon(u_\varepsilon) = 0$, we denote the corresponding Esscher transform of $\bar{\mathbb{P}}^{\varepsilon}$ by $\bar{\mathbb{Q}}^{\varepsilon}$ and conclude, using Proposition 1.3, that

$$\mathbb{P}(\|X\| \leq 3\varepsilon) \geq \mathbb{P}(X \text{ has no jumps larger than } \varepsilon) \cdot \bar{\mathbb{P}}^{\varepsilon}(\|X\| \leq 3\varepsilon)$$
$$\geq \exp\{-\nu([-\varepsilon,\varepsilon]^c) + \Lambda_\varepsilon(u_\varepsilon) - 3\varepsilon|u_\varepsilon|\} \cdot \bar{\mathbb{Q}}^{\varepsilon}(\|X\| \leq 3\varepsilon).$$

Under $\bar{\mathbb{Q}}^{\varepsilon}$, the process $X$ is a $(e^{u_\varepsilon x} \cdot \nu(dx)|_{[-\varepsilon,\varepsilon]}, \sigma^2)$-Lévy martingale. Denoting by $\bar{F}$ the corresponding $F$-function from Proposition 1.4, that is,

$$(14) \qquad \bar{F}(\varepsilon) = \frac{1}{\varepsilon^2}\left[\sigma^2 + \int_{-\varepsilon}^{\varepsilon} x^2 e^{u_\varepsilon x}\nu(dx)\right] = \varepsilon^{-2}\Lambda''_\varepsilon(u_\varepsilon),$$

we conclude with Proposition 1.4 that

$$\mathbb{P}(\|X\| \leq 3\varepsilon) \geq \exp\{-[\nu([-\varepsilon,\varepsilon]^c) - \Lambda_\varepsilon(u_\varepsilon) + 3\varepsilon|u_\varepsilon| + 10\bar{F}(\varepsilon) + 3]\}.$$

Conversely, analog computations show that

$$\mathbb{P}(\|X\| \leq \varepsilon/2) \leq \exp\left\{-\left[\nu([-\varepsilon,\varepsilon]^c) - \Lambda_\varepsilon(u_\varepsilon) - \frac{\varepsilon}{2}|u_\varepsilon| + \frac{1}{12}\bar{F}(\varepsilon) - 1\right]\right\}.$$

Our main theorem summarizes these considerations:

THEOREM 1.5. *Let $\varepsilon > 0$ and let $X$ be a Lévy process with triplet $(\nu, \sigma^2, b)$ that possesses a solution $u_\varepsilon \in \mathbb{R}$ to the equation $\Lambda'_\varepsilon(u_\varepsilon) = 0$, where $\Lambda_\varepsilon$ is as in (13). Then, one has*

$$(15) \qquad -\log\mathbb{P}(\|X\| \leq 3\varepsilon) \leq \nu([-\varepsilon,\varepsilon]^c) - \Lambda_\varepsilon(u_\varepsilon) + 3\varepsilon|u_\varepsilon| + 10\bar{F}(\varepsilon) + 3$$

*and*

$$(16) \qquad -\log \mathbb{P}(\|X\| \leq \varepsilon/2) \geq \nu([-\varepsilon,\varepsilon]^c) - \Lambda_\varepsilon(u_\varepsilon) - \frac{\varepsilon}{2}|u_\varepsilon| + \frac{1}{12}\bar{F}(\varepsilon) - 1,$$

*where* $\bar{F}(\varepsilon)$ *is as in (14).*

A particularly important consequence of Theorem 1.5 is that for Lévy processes with nonvanishing Gaussian component $\sigma$, one has

$$-\log \mathbb{P}(\|X\| \leq \varepsilon) \sim \frac{\pi^2}{8}\sigma^2 \varepsilon^{-2};$$

see Corollary 2.6.

Theorem 1.5 gives estimates of (1) in terms of $\nu([-\varepsilon,\varepsilon]^c)$, $\Lambda_\varepsilon(u_\varepsilon)$, $|u_\varepsilon|$ and $\bar{F}(\varepsilon)$. These quantities depend, in a nontrivial way, on the characterizing triplet of the Lévy process. We want to emphasize that the lower and upper bound are tight in the sense that

$$(17) \qquad -\log \mathbb{P}(\|X\| \leq 3\varepsilon) \preceq \nu([-\varepsilon,\varepsilon]^c) - \Lambda_\varepsilon(u_\varepsilon) + \bar{F}(\varepsilon)$$

$$\preceq -\log \mathbb{P}(\|X\| \leq \varepsilon/2),$$

whenever

$$(18) \qquad \varepsilon|u_\varepsilon| \prec -\Lambda_\varepsilon(u_\varepsilon) + \nu([-\varepsilon,\varepsilon]^c) + \bar{F}(\varepsilon),$$

where $f \prec g$ means $\limsup f/g = 0$. This condition is satisfied in all examples considered below and we are not aware of any counterexamples. Furthermore, we give sufficient conditions for (18) to hold in the Appendix; see Lemma A.1.

Note that equation (17) gives the weak asymptotic order of the small deviations (1) whenever $-\log \mathbb{P}(\|X\| \leq 2\varepsilon) \approx -\log \mathbb{P}(\|X\| \leq \varepsilon)$ and condition (18) is satisfied.

Interestingly, for general Lévy processes, the probability of a small deviation can be arbitrarily small [thus, the expression in (1) can increase arbitrarily fast; see Remark 2.4] so that one can easily construct examples for which $-\log \mathbb{P}(\|X\| \leq 2\varepsilon) \approx -\log \mathbb{P}(\|X\| \leq \varepsilon)$ fails to hold.

Assuming the validity of (18), we see that the small deviations are governed by three effects:

- the first term, that is, $\nu([-\varepsilon,\varepsilon]^c)$, represents the cost of having *no large jumps*;
- the second term, that is, $-\Lambda_\varepsilon(u_\varepsilon)$, represents the cost induced by the *drift* of the modified process (the *Esscher term*);
- the third term, that is, $\bar{F}(\varepsilon)$, represents the cost induced by the *oscillations* of the modified process.

We remark that any of the terms $\nu([-\varepsilon,\varepsilon]^c)$, $-\Lambda_\varepsilon(u_\varepsilon)$ and $\bar{F}(\varepsilon)$ can give the leading term in the asymptotics. An explicit analysis is carried out below for many examples.

**2. Examples.** The steps described in the last subsection enable us to derive an estimate for the rate in (1) for any given Lévy process. Let us demonstrate the technique with some examples. In particular, we show how to re-prove all known results about the small deviation order for Lévy processes using our framework.

2.1. *Symmetric Lévy processes.* The first example concerns the case when the distribution of $X$ is symmetric. Then, $X$ is already a $\bar{\mathbb{P}}^\varepsilon$-martingale and we get the following simpler bounds:

COROLLARY 2.1. *Let $X$ be a symmetric Lévy process [i.e., $\mathcal{L}(X_1) = \mathcal{L}(-X_1)$]. Then,*

$$-\log \mathbb{P}(\|X\| \leq \varepsilon) \approx \nu([-\varepsilon, \varepsilon]^c) + F(\varepsilon),$$

*where $F$ is as defined in (12).*

The corollary follows immediately from Theorem 1.5, once one notices that, for any Lévy measure $\nu$ and Gaussian component,

$$\nu([-2\varepsilon, 2\varepsilon]^c) + F(2\varepsilon) \approx \nu([-\varepsilon, \varepsilon]^c) + F(\varepsilon).$$

We give a proof of this fact in Lemma 5.1. Let us concretize the last result when the Lévy measure is given by some regularly varying function.

EXAMPLE 2.2. *Let $X$ be a symmetric Lévy process with $\sigma^2 = 0$ and Lévy measure*

$$\nu([-\varepsilon, \varepsilon]^c) \approx \varepsilon^{-\alpha} \ell(\varepsilon) \qquad as \ \varepsilon \to 0$$

*for some slowly varying function $\ell$ and $0 < \alpha \leq 2$ [note that, for $\alpha = 2$, certain restrictions for $\ell$ apply in order to ensure (3)]. It is then easily seen that*

$$F(\varepsilon) \approx \varepsilon^{-2} \int_0^\varepsilon x^{1-\alpha} \ell(x) \, dx.$$

*If $\alpha < 2$ and $\ell$ is a slowly varying function that is bounded away from 0 and $\infty$ on any compact interval, then the last term behaves asymptotically as $\varepsilon^{-\alpha} \ell(\varepsilon)$; see [4]. However, this is not true for $\alpha = 2$. Namely, let us evaluate $F$ when $\ell(x) = c |\log x|^{-\gamma}$. We then obtain*

$$-\log \mathbb{P}(\|X\| \leq \varepsilon) \approx \begin{cases} \varepsilon^{-\alpha} |\log \varepsilon|^{-\gamma}, & 0 < \alpha < 2, \\ \varepsilon^{-2} |\log \varepsilon|^{1-\gamma}, & \alpha = 2, \gamma > 1. \end{cases}$$

*The above example includes Lévy processes that are approximately (in the sense that the asymptotic behavior of the Lévy measure at zero is the same) a symmetric $\alpha$-stable Lévy process.*

2.2. *Subordinators with negative drift.* We now consider a class of subordinators with additional negative drift.

EXAMPLE 2.3. *Let us first return to the stable subordinator $X$ with $0 < \alpha < 1$ considered in (8), where we add the drift with slope $\mu$. The cases $\mu > 0$ and $\mu = 0$ were treated in Example 1.2, so assume that $\mu < 0$. In this case, Theorem 1.5 yields the somewhat surprising result [recall (6) for the subordinator without drift and (5) for the strictly stable case]*

$$(19) \qquad -\log \mathbb{P}\Big( \sup_{t \in [0,1]} |X_t + \mu t| \le \varepsilon \Big) \approx \varepsilon^{-1}|\log \varepsilon|.$$

The following remark illustrates that expressions of the form $\varepsilon^{-1}|\log \varepsilon|$ appear naturally in the context of small deviations due to their relation to the *large deviations* of the Poisson distribution.

REMARK 2.4. Assume, now, that $\nu$ is a finite measure that is supported on $(0, \infty)$ and set $f(\varepsilon) = \nu(0, \varepsilon]$ for $\varepsilon > 0$. We consider a $(\nu, 0, -1 + \int_0^1 x\nu(dx))$-Lévy process $X$. As before, we estimate

$$\mathbb{P}(\|X\| \le \varepsilon/2) \le \mathbb{P}(X \text{ has no jumps larger than } \varepsilon) \cdot \bar{\mathbb{P}}^\varepsilon(X_1 \ge -\varepsilon).$$

Let $N_\varepsilon$ denote the number of jumps smaller or equal to $\varepsilon$. Then, $X_1 \le -1 + \varepsilon N_\varepsilon$, $\bar{\mathbb{P}}^\varepsilon$-almost surely, so that

$$\bar{\mathbb{P}}^\varepsilon(X_1 \ge -\varepsilon) \le \mathbb{P}\Big( N_\varepsilon \ge \frac{1}{\varepsilon} - 1 \Big).$$

The random variable $N_\varepsilon$ is Poisson distributed with parameter $f(\varepsilon)$ and we derive, by using the exponential Chebyshev inequality, that

$$
\begin{aligned}
(20) \qquad \mathbb{P}(\|X\| \le \varepsilon/2) &\le \mathbb{P}\Big( N_\varepsilon \ge \frac{1}{\varepsilon} - 1 \Big) \\
&\le \exp\Big( -\Big[\Big(\frac{1}{\varepsilon} - 1\Big)\Big(\log \frac{1/\varepsilon - 1}{f(\varepsilon)} - 1\Big) + f(\varepsilon)\Big]\Big) \\
&= \exp\Big( -(1 + o(1))\frac{1}{\varepsilon}[|\log \varepsilon| + |\log f(\varepsilon)|]\Big).
\end{aligned}
$$

Thus, we always get that

$$-\log \mathbb{P}(\|X\| \le \varepsilon) \succeq \varepsilon^{-1}|\log \varepsilon|.$$

Moreover, the estimate (20) shows that one can achieve arbitrarily small probabilities [i.e., arbitrarily fast increase of (1)] for the small deviation by choosing $\nu$ accordingly.

EXAMPLE 2.5. Our next example concerns the Gamma process $X$:

$$\mathbb{E}e^{izX_t} = \exp\left(t\int_0^\infty (e^{izx}-1)\frac{e^{-x/a}b\,dx}{x}\right),$$

with parameters $a, b > 0$. We add a drift with slope $\mu$. For $\mu > 0$, the process $(X_t + \mu t)$ does not satisfy the small deviation property. For $\mu = 0$, we are in the case of a subordinator with $\sup_{t\in[0,1]} |X_t| = X_1$ being Gamma-distributed so that

$$\mathbb{P}\left(\sup_{t\in[0,1]} |X_t| \le \varepsilon\right) \sim \frac{\varepsilon^{b+1}}{a^b\Gamma(b+1)}.$$

On the other hand, for a negative drift $\mu < 0$, our results imply that

$$-\log\mathbb{P}\left(\sup_{t\in[0,1]} |X_t + \mu t| \le \varepsilon\right) \approx \varepsilon^{-1}|\log\varepsilon|.$$

2.3. *Lévy processes with nonvanishing Gaussian component.* Let us now look at Lévy processes with nonvanishing Gaussian component (i.e., assume that $\sigma \ne 0$).

COROLLARY 2.6. *Let $X$ be a $(\nu, \sigma, b)$-Lévy process with $\sigma \ne 0$. Then,*

$$-\log\mathbb{P}(\|X\| \le \varepsilon) \sim \frac{\pi^2}{8}\sigma^2\varepsilon^{-2}.$$

PROOF. We represent $X$ as $X = Y + \sigma B$, where $B$ is a standard Brownian motion, $Y$ does not contain any Gaussian component, and $B$ and $Y$ are independent. By Anderson's inequality [2],

$$\mathbb{P}(\|Y + \sigma B\| \le \varepsilon) \le \mathbb{P}(\|\sigma B\| \le \varepsilon),$$

which already implies, by (4), that

$$\liminf_{\varepsilon\to 0} \varepsilon^2(-\log\mathbb{P}(\|Y + \sigma B\| \le \varepsilon)) \ge \frac{\pi^2}{8}\sigma^2.$$

On the other hand, let $0 < \varepsilon < 1$, $0 < \theta < 1$ and let

$$X_t = Y_t - b_\varepsilon t + \sigma B_t + b_\varepsilon t + Z_t, \qquad t \in [0,1],$$

where $\sigma B$ is the Gaussian component of $X$, $Z$ is constructed from the jumps of $X$ that are larger than $\varepsilon$ and $b_\varepsilon$ is chosen such that $Y_t - b_\varepsilon t$ is a martingale, that is,

$$b_\varepsilon := b - \int_{\{\varepsilon < |x| < 1\}} x\nu(dx).$$

Observe that the components $Y$, $B$ and $Z$ are independent and, thus,

$$
\begin{aligned}
(21) \quad & \mathbb{P}(\|Y_t - b_\varepsilon t\| \le \theta \varepsilon) \mathbb{P}(\|\sigma B_t + b_\varepsilon t\| \le (1-\theta)\varepsilon) \mathbb{P}(\|Z\| = 0) \\
& = \mathbb{P}(\|Y_t - b_\varepsilon t\| \le \theta \varepsilon, \|\sigma B_t + b_\varepsilon t\| \le (1-\theta)\varepsilon, \|Z\| = 0) \\
& \le \mathbb{P}(\|X\| \le \varepsilon).
\end{aligned}
$$

Clearly,

$$
-\log \mathbb{P}(\|Z\| = 0) = \nu(\{|x| > \varepsilon\}) = o(\varepsilon^{-2}).
$$

Furthermore, by Proposition 1.4,

$$
-\log \mathbb{P}(\|Y_t - b_\varepsilon t\| \le \theta \varepsilon) \approx \varepsilon^{-2} \int_{\{|x| \le \theta \varepsilon\}} x^2 \nu(dx) = o(\varepsilon^{-2}).
$$

On the other hand, by Proposition 1.3 and (4),

$$
-\log \mathbb{P}(\|\sigma B_t + b_\varepsilon t\| \le (1-\theta)\varepsilon) \lesssim \frac{b_\varepsilon^2}{2\sigma^2} - \log \mathbb{P}(\|\sigma B\| \le (1-\theta)\varepsilon)
$$

$$
\lesssim (1-\theta)^{-2} \frac{\pi^2}{8} \sigma^2 \varepsilon^{-2}
$$

since $|b_\varepsilon| = o(\varepsilon^{-1})$. Therefore, (22) implies that

$$
(22) \qquad \limsup_{\varepsilon \to 0} \varepsilon^2 (-\log \mathbb{P}(\|X\| \le \varepsilon)) \le (1-\theta)^{-2} \frac{\pi^2}{8} \sigma^2.
$$

Letting $\theta$ tend to zero completes the proof. $\quad \square$

2.4. *Polynomial Lévy measure.* Let us look at what happens if both lower tails of the Lévy measure are polynomial with different exponents. The technique used for this example can be extended to any case with regularly varying Lévy measure at zero.

Let $X$ be a Lévy process with triplet $(\nu, 0, b)$, where $\nu = \nu_0 + \nu_1$, $\nu_1$ is some finite measure concentrated on $\{|x| > 1\}$ and $\nu_0$ is given by

$$
(23) \qquad \frac{\nu_0(dx)}{dx} = \frac{C_1 \mathbb{1}_{(0,1]}(x)}{x^{1+\alpha_1}} + \frac{C_2 \mathbb{1}_{[-1,0)}(x)}{(-x)^{1+\alpha_2}},
$$

where $\alpha_1, \alpha_2 < 2$ and $C_1, C_2 \ge 0$, $C_1 + C_2 \ne 0$.

First, let us note that $\nu_1$ has no influence on the order, so we can, and will, assume without loss of generality that $\nu_1 = 0$. On the other hand, observe that, if $\alpha_1 \ne \alpha_2$, we can always assume $\alpha_1 > \alpha_2$ (by passing over to $-X$ if necessary). Equally, if $\alpha_1 = \alpha_2$, we can assume that $C_1 > C_2$, unless we are in the case of symmetric $\nu$. This reduces the number of cases that have to be treated.

We distinguish three regimes: the cases where $\alpha_1 > 1$, $\alpha_1 = 1$ and $0 < \alpha_1 < 1$. The second exponent $\alpha_2$ can even be negative.

COROLLARY 2.7.    *Let $\alpha_1 \geq \alpha_2$ and $\alpha_1 > 1$. Then,*

$$-\log \mathbb{P}(\|X\| \leq \varepsilon) \approx \varepsilon^{-\alpha_1}.$$

If the maximal exponent equals 1, then we get the following:

COROLLARY 2.8.    *Let $\alpha_1 > \alpha_2$ and $\alpha_1 = 1$. Then,*

$$-\log \mathbb{P}(\|X\| \leq \varepsilon) \approx \varepsilon^{-1} |\log \varepsilon| \log(|\log \varepsilon|).$$

*The same holds in the case $\alpha_1 = \alpha_2 = 1$ and $C_1 > C_2$. However, the case $\alpha_1 = \alpha_2 = 1$ and $C_1 = C_2$ leads to*

$$-\log \mathbb{P}(\|X\| \leq \varepsilon) \approx \varepsilon^{-1}.$$

Finally, if the maximal exponent is less than 1, we obtain the following result:

COROLLARY 2.9.    *Let $\alpha_1 \geq \alpha_2$ and $0 < \alpha_1 < 1$. Set $c := b - \int_{-1}^{1} x\nu(dx)$. Then,*

$$-\log \mathbb{P}(\|X\| \leq \varepsilon) \approx \begin{cases} \varepsilon^{-1} |\log \varepsilon|, & c \neq 0, \\ \varepsilon^{-\alpha_1}, & c = 0. \end{cases}$$

Let us briefly discuss our findings.

REMARK 2.10.    In the first regime, the case $\alpha_1 > 1$, the asymmetry does not have any influence on the small deviations. In particular, the cost for having no large jumps and the cost induced by the oscillations govern the asymptotics.

In the case $\alpha_1 = 1$, the result strongly depends on the magnitude of the asymmetry. If $\alpha_1$ is not equal to $\alpha_2$ or $C_1$ is not equal to $C_2$, then the Esscher term dominates the small deviations. Otherwise, we regain the asymptotics of the symmetric case. Note that, so far, the drift $b$ has not influenced the asymptotics.

In the case $\alpha_1 < 1$, we have $\int_{-1}^{1} |x|\nu(dx) < \infty$; thus, we can define the *effective drift*, $c = b - \int_{-1}^{1} x\nu(dx)$. The result now strongly depends on this effective drift. If this is nonzero, then the Esscher term determines the small deviation order. Otherwise, the nonexistence of large jumps and the oscillation term govern the asymptotics.

The case of symmetric Lévy measure (i.e., $\alpha_1 = \alpha_2 =: \alpha$ and $C_1 = C_2$) is included in the above results; however, we repeat it here since it represents the most important situation.

COROLLARY 2.11. *Let $\alpha_1 = \alpha_2 =: \alpha$, $C_1 = C_2$. If $b \neq 0$, then*

$$-\log \mathbb{P}(\|X\| \leq \varepsilon) \approx \begin{cases} \varepsilon^{-1} |\log \varepsilon|, & 0 < \alpha < 1, \\ \varepsilon^{-\alpha}, & 1 \leq \alpha < 2. \end{cases}$$

*If $b = 0$, we are in the symmetric case and the rate is $\varepsilon^{-\alpha}$ for all $\alpha \in (0, 2)$.*

The latter corollaries cover, in particular, general $\alpha$-stable Lévy processes, that is, if $\alpha_1 = \alpha_2 =: \alpha$ and $\nu_1$ is adjusted appropriately, we deal with an $\alpha$-stable Lévy process. This generalizes the known results for strictly stable processes (5) and stable subordinators (6). Also, tempered $\alpha$-stable processes (see, e.g., [17]) are included in the above results.

### 2.5. *Further examples.*

EXAMPLE 2.12. Let us consider the so-called variance Gamma process. This process is obtained when replacing the time parameter of a Brownian motion with drift by a Gamma subordinator, that is, letting $B$ be a Brownian motion, $\sigma > 0$, $\mu \in \mathbb{R}$, and $A$ be a Gamma process as defined in Example 2.5, independent of $B$. Then, $X_t := \sigma B_{A_t} + \mu A_t$ is called a variance Gamma process. It is a Lévy process with Lévy measure $\nu$ given by

$$\frac{\nu(dx)}{dx} = \frac{C_1}{x} e^{-\lambda_1 x} \mathbb{1}_{\{x > 0\}} + \frac{C_2}{(-x)} e^{-\lambda_2(-x)} \mathbb{1}_{\{x < 0\}}$$

with parameters $C_1, C_2, \lambda_1, \lambda_2 > 0$, depending in some way on $\sigma$, $\mu$ and the parameters of the Gamma process. In particular, $C_1 = C_2$ and $\lambda_1 = \lambda_2$ if and only if $\mu = 0$.

Applying Theorem 1.5, we see that the small deviation probability of $X$ is given by

$$-\log \mathbb{P}(\|X\| \leq \varepsilon) \approx \begin{cases} |\log \varepsilon|, & \mu = 0, \\ \varepsilon^{-1} |\log \varepsilon|, & \mu \neq 0. \end{cases}$$

This extends the recent results for the case $\mu = 0$ in [15]. When $\mu \neq 0$, the Esscher term dominates. In the case $\mu = 0$, the very slow increase comes from the very slow increase of the Gamma process; see [15].

EXAMPLE 2.13. In [14], the small deviation problem was solved for processes arising from subordination to fractional Brownian motion. For the case of subordination to Brownian motion, the resulting process is a Lévy process and the findings follow easily from our current framework.

In order to formulate the result, let $B$ be a Brownian motion and $A$ be any subordinator (independent of $B$) with Laplace exponent $\Phi$, that is,

$$\Phi(u) := -\log \mathbb{E} e^{-uA_1} = u b_A + \int_0^\infty (1 - e^{-ux}) \nu_A(dx),$$

where $\nu_A$ is the Lévy measure of the subordintor $A$ and $b_A \geq 0$ is the drift of $A$. We then consider the subordinated process $X_t := B_{A_t}$. Note that $X$ is symmetric, so, by Corollary 2.1, we only have to investigate $\nu([-\varepsilon, \varepsilon]^c)$ and $F$ from (12), where $\nu$ is the Lévy measure of $X$. Doing so, we find that

$$-\log \mathbb{P}(\|X\| \leq \varepsilon) \approx \Phi(\varepsilon^{-2}) + b_A \left\{ \varepsilon^{-2} \left( \sigma^2 + \int_0^\varepsilon x^2 \nu_A(dx) \right) + \nu_A([-\varepsilon, \varepsilon]^c) \right\}$$

for any subordinator $A$. This improves the results from [14] since, there, certain regularity conditions for the Laplace exponent $\Phi$ and $b_A = 0$ were assumed. Again, the case of subordination with the Gamma process ($\Phi(u) = b \log(u + 1/a)$, $b_A = 0$, treated in [15]) is included.

EXAMPLE 2.14. Let us now consider a compound Poisson process with no effective drift (for some remarks on compound Poisson processes with drift, see Remark 2.4), that is, Lévy processes with finite Lévy measure $\nu$ and $c := b - \int_{-1}^1 x \nu(dx) = 0$. Here, the small deviation probability does not tend to zero [i.e., (1) does not tend to $\infty$] because the probability that the compound Poisson process has no jump is positive. In fact, we have $\mathbb{P}(\sup_{t \in [0,1]} |X_t| \leq \varepsilon) \to \exp(-\nu(\mathbb{R}))$.

**3. Esscher term.** In this section, we prove Proposition 1.3. We assume that $X$ is a Lévy process with compactly supported Lévy measure $\nu$. As outlined above, this can be obtained by subtracting the large jumps from the Lévy process.

The goal of this section is to show that one can transform the Lévy process by using the Esscher transform in such a way that the resulting process is a Lévy martingale. This transformation incurs a "cost" on the probability that the process remains inside an $\varepsilon$-strip. The resulting Lévy martingale can then be treated with the methods in Section 4.

Furthermore, we prove in this section that the Esscher transformation in Proposition 1.3 is possible in all cases of interest (Lemma 3.2). However, first, we prove Proposition 1.3.

PROOF OF PROPOSITION 1.3. We now denote by $X$ a Lévy process with compactly supported Lévy measure $\nu$ (we exclude the trivial case when $\nu = 0$, $b = 0$ and $\sigma^2 = 0$) and represent its characteristic function as $\mathbb{E} e^{iuX_t} = e^{-t\psi(u)}$, where

$$\psi(u) = \frac{\sigma^2}{2} u^2 - ibu + \int_{\mathbb{R} \setminus \{0\}} (1 - e^{iux} + iux) \nu(dx).$$

For $u \in \mathbb{R}$, we consider the Esscher transform $\mathbb{Q}^u$ given by

$$\left. \frac{d\mathbb{Q}^u}{d\mathbb{P}} \right|_{\mathcal{F}_t} = e^{uX_t - t\Lambda(u)}.$$

Here, $(\mathcal{F}_t)$ denotes the canonical filtration induced by the process $X$. Observe that

$$\log \mathbb{E}^{\mathbb{Q}^u} e^{\theta X_t} = \log \mathbb{E} e^{(\theta+u)X_t - t\Lambda(u)} = t(\Lambda(\theta+u) - \Lambda(u)).$$

Using the fact that $\Lambda(u) = \frac{1}{2}\sigma^2 u^2 + bu + \int [e^{ux} - 1 - ux]\nu(dx)$, we conclude that

$$
\begin{aligned}
\Lambda(\theta+u) - \Lambda(u) &= \frac{1}{2}\sigma^2(\theta+u)^2 + b(\theta+u) \\
&\quad + \int [e^{(\theta+u)x} - 1 - (\theta+u)x]\nu(dx) \\
&\quad - \left[\frac{1}{2}\sigma^2 u^2 + bu + \int [e^{ux} - 1 - ux]\nu(dx)\right] \\
&= \frac{1}{2}\sigma^2\theta^2 + \left(b + \sigma^2 u + \int x(e^{ux} - 1)\nu(dx)\right)\theta \\
&\quad + \int [e^{\theta x} - 1 - \theta x]e^{ux}\nu(dx).
\end{aligned}
$$

Thus, $X$ is a $\mathbb{Q}^u$-Lévy process. Moreover, it is an $(e^{u^*x} \cdot \nu, \sigma^2)$-Lévy martingale if $u$ is equal to the solution $u^*$ of the equation

$$\Lambda'(u^*) = b + \sigma^2 u^* + \int x(e^{u^*x} - 1)\nu(dx) = 0.$$

In the sequel, we let $\mathbb{Q} := \mathbb{Q}^{u^*}$. Then,

$$e^{\Lambda(u^*) - \varepsilon|u^*|}\mathbb{Q}(\|X\| \le \varepsilon) \le \mathbb{P}(\|X\| \le \varepsilon) \le e^{\Lambda(u^*) + \varepsilon|u^*|}\mathbb{Q}(\|X\| \le \varepsilon). \qquad \square$$

REMARK 3.1. The property $\Lambda'(u^*) = 0$, together with the convexity of $\Lambda$, implies that

$$\Lambda(u^*) = \inf_{u \in \mathbb{R}} \Lambda(u).$$

Thus, the change of measure leads to an equivalent martingale measure that is entropy-minimizing.

Proposition 1.3 implies that the existence of a solution to $\Lambda'(u) = 0$ yields the existence of a so-called equivalent martingale measure $\mathbb{Q}$ for $X$. Certainly, such equivalent martingale measures do not always exist. In particular, all subordinators do not possess equivalent martingale measures and thus we will not be able to apply Proposition 1.3 for these processes.

Fortunately, the case of $X$ being a subordinator is essentially the only relevant case in which we cannot apply Proposition 1.3:

LEMMA 3.2. *Let $X$ be a Lévy process. If $|X|$ is not a subordinator and possesses the small deviation property, then (9) has a unique solution under each measure $\overline{\mathbb{P}}^\varepsilon$.*

The last lemma shows that either:

(a) the problem (9) has a solution with $u^* \in \mathbb{R}$—and we can thus work with the estimate (11); or

(b) the process is a subordinator without drift (or the negative of a subordinator)—in which case the small deviation problem is solved via the Tauberian theorem; or

(c) the process does not satisfy the small deviation property—which means that for sufficiently small $\varepsilon$, the probability is zero.

PROOF OF LEMMA 3.2. Clearly, we can assume that $0 < \varepsilon < 1$ and $\nu(|x| > \varepsilon) = 0$. We have to check whether the function

$$\Lambda'(u) = \sigma^2 u + b + \int_{-\varepsilon}^{\varepsilon} (e^{ux} - 1) x \nu(dx)$$

has a unique root. First, we note that $u \mapsto (e^{ux} - 1)x$ is a strictly increasing function for any $x \in \mathbb{R} \setminus \{0\}$. Thus, the function $u \mapsto \int_{-\varepsilon}^{\varepsilon} (e^{ux} - 1)x\nu(dx)$ is strictly increasing (unless $\nu = 0$, which is trivial) and is hence $\Lambda'$. This makes any root unique. Furthermore, we note that $\Lambda'$ is continuous so that existence of a root is equivalent to $\lim_{u \to -\infty} \Lambda'(u) < 0 < \lim_{u \to \infty} \Lambda'(u)$.

Let us consider the various cases.

*Case* 1: $\sigma \neq 0$. In this case, we clearly have $\lim_{u \to \pm\infty} \Lambda'(u) = \pm\infty$.

*Case* 2: We consider $\sigma = 0$ and $\int_{-\varepsilon}^{\varepsilon} |x| \nu(dx) = \infty$. Then, we must have $\int_0^\varepsilon x\nu(dx) = \infty$ or $\int_{-\varepsilon}^0 -x\nu(dx) = \infty$. In the former case,

$$\lim_{u \to \pm\infty} \int_0^\varepsilon (e^{ux} - 1)x\nu(dx) \to \pm\infty \quad \text{and}$$

$$\int_{-\varepsilon}^0 (e^{ux} - 1)x\nu(dx) \begin{cases} \geq 0, & u > 0, \\ \leq 0, & u < 0, \end{cases}$$

which shows the assertion. The latter case is treated analogously.

*Case* 3: Finally, let $X$ be of type (I). This is the most interesting case, where we actually need the assumptions. In this case, we have

$$\Lambda'(u) = b - \int_{-\varepsilon}^{\varepsilon} x\nu(dx) + \int_0^\varepsilon e^{ux} x\nu(dx) + \int_{-\varepsilon}^0 e^{ux} x\nu(dx).$$

*Case* 3(a): If $\nu(0 < x < \varepsilon) > 0$ and $\nu(-\varepsilon < x < 0) > 0$, then clearly $\lim_{u \to \pm\infty} \Lambda'(u) = \pm\infty$.

*Case* 3(b): Let $\nu(0 < x < \varepsilon) > 0$ and $\nu(-\varepsilon < x < 0) = 0$. Then, $\lim_{u\to\infty} \Lambda'(u) = \infty$ and $\lim_{u\to-\infty} \Lambda'(u) = b - \int_{-\varepsilon}^{\varepsilon} x\nu(dx) = -b'$. If the latter term is positive, then the process does not satisfy the small deviation property (by Proposition 1.1) and we are done. If $-b'$ equals zero, then it is easily seen that the process is in fact a subordinator. And, finally, if it is negative, then $\Lambda'$ must have a (unique) root.

Case 3(c), where $\nu(0 < x < \varepsilon) = 0$ and $\nu(-\varepsilon < x < 0) > 0$, is treated as 3(b). Case 3(d), where $\nu(0 < x < \varepsilon) = 0$ and $\nu(-\varepsilon < x < 0) = 0$, is trivial. $\square$

**4. Exit time arguments.** In this section, we consider a Lévy martingale having jumps smaller than $\varepsilon$. Note that one can obtain this by first removing the jumps larger than $\varepsilon$ and then applying the transformation described in Section 3.

The purpose of this section is to prove Proposition 1.4. We proceed in several steps. First, let us define inductively the times

$$\tau_i = \inf\{t \geq 0 : |X_{\tau_1 + \cdots + \tau_{i-1} + t} - X_{\tau_1 + \cdots + \tau_{i-1}}| \geq \varepsilon\}, \qquad i = 1, 2, \ldots,$$

and the increments $Z_i = X_{\tau_1 + \cdots + \tau_i} - X_{\tau_1 + \cdots + \tau_{i-1}}$. To be formally correct, we here need to consider the Lévy process defined on the whole interval $[0, \infty)$. By the strong Markov property, the family $(Z_i, \tau_i)_{i\in\mathbb{N}}$ consists of independent identically distributed random variables.

LEMMA 4.1. *For $\varepsilon > 0$, we have*

$$\mathbb{P}(\|X\| \leq 3\varepsilon) \geq e^{-10F(\varepsilon)-3}.$$

PROOF. For convenience, we let $\tau = \tau_1$ and $Z = Z_1$. Due to the optional stopping theorem and the boundedness of the jumps, we have

$$0 = \mathbb{E}X_\tau \leq 2\varepsilon\mathbb{P}(Z > 0) - \varepsilon\mathbb{P}(Z < 0)$$

so that

$$\mathbb{P}(Z > 0) \geq 1/3.$$

Moreover, Doob's martingale inequality gives that

$$\mathbb{P}(\tau \leq t) = \mathbb{P}\Big(\sup_{s\leq t}|X_s| \geq \varepsilon\Big) \leq \frac{\mathbb{E}|X_t|^2}{\varepsilon^2} = F(\varepsilon)t.$$

Hence, we have, for $t := (4F(\varepsilon))^{-1}$,

$$\mathbb{P}(\tau \geq t, Z > 0) \geq 1 - \mathbb{P}(\tau \leq t) - \mathbb{P}(Z < 0) \geq \tfrac{1}{12}.$$

By symmetry, we also have $\mathbb{P}(\tau \geq t, Z < 0) \geq 1/12$.

Next, fix the smallest integer $n$ with $n \geq t^{-1}$ and consider the event

$$E = \{\forall i \in \{1, \ldots, n\} : \tau_i \geq t, \operatorname{sgn}(Z_i) = -\operatorname{sgn}(X_{\tau_1 + \cdots + \tau_i - 1})\},$$

where sgn denotes the signum function, with $\operatorname{sgn}(0) = 1$. If $E$ occurs, then $(X_t)$ starts at each time $\tau_1 + \cdots + \tau_{i-1}$ in the interval $[-2\varepsilon, 2\varepsilon]$ and ends at time $\tau_1 + \cdots + \tau_i$ in the same interval. Hence, along the whole trajectory, we have $|X_s| \leq 3\varepsilon$ while $s \leq \sum_{i=1}^{n} \tau_i$. Since, by assumption, $\sum_{i=1}^{n} \tau_i \geq nt \geq 1$, we have $\|X\| \leq 3\varepsilon$. Hence, using the strong Markov property of the Lévy process, we obtain

$$\mathbb{P}(\|X\| \leq 3\varepsilon) \geq \mathbb{P}(E) \geq 12^{-n} \geq 12^{-4F(\varepsilon)-1} \geq e^{-10F(\varepsilon)-3}. \qquad \square$$

Conversely, one can prove the following lemma.

LEMMA 4.2.  *For $\varepsilon > 0$, we have*

$$\mathbb{P}(\|X\| \leq \varepsilon/2) \leq e^{-F(\varepsilon)/12+1}.$$

PROOF.  Again, let $\tau$ denote the first exit time of $X$ out of $[-\varepsilon, \varepsilon]$. Then, by Wald's identity,

$$4\varepsilon^2 \geq \limsup_{t \to \infty} \mathbb{E} X_{t \wedge \tau}^2 = \limsup_{t \to \infty} \varepsilon^2 F(\varepsilon) \mathbb{E}[t \wedge \tau] = \varepsilon^2 F(\varepsilon) \mathbb{E}\tau$$

and thus, by the Markov inequality,

$$\mathbb{P}(\tau \geq 8/F(\varepsilon)) \leq 1/2.$$

Consequently, one has, for $n := \lfloor F(\varepsilon)/8 \rfloor > 0$ and $t_i = 8i/F(\varepsilon)$, $i = 0, \ldots, n$,

$$\mathbb{P}(\|X\| \leq \varepsilon/2) \leq \mathbb{P}\Big(\forall i = 0, \ldots, n-1 : \sup_{s \in [t_i, t_{i+1}]} |X_s - X_{t_i}| \leq \varepsilon\Big)$$

$$\leq 2^{-n} \leq 2^{-F(\varepsilon)/8+1}.$$

If $n = 0$, the result holds trivially.   $\square$

**5. Proofs of the explicit rates in the examples.**  In this section, we give the proofs for the asymptotic rates stated in the examples.

First, we show the following lemma, which immediately yields Corollary 2.1.

LEMMA 5.1.  *With $F$ as in (12), we have*

$$\nu([-2\varepsilon, 2\varepsilon]^c) + F(2\varepsilon) \leq \nu([-\varepsilon, \varepsilon]^c) + F(\varepsilon) \leq 4[\nu([-2\varepsilon, 2\varepsilon]^c) + F(2\varepsilon)].$$

PROOF. For $\varepsilon > 0$, we consider the function $g_\varepsilon : [0, \infty) \to [0, 1]$ defined by $g_\varepsilon(x) = \frac{x^2}{\varepsilon^2} \wedge 1$. Then, $g_{2\varepsilon} \leq g_\varepsilon \leq 4g_{2\varepsilon}$ so that

$$\nu([-\varepsilon, \varepsilon]^c) + F(\varepsilon) = \int g_\varepsilon \, d\nu + \frac{\sigma^2}{\varepsilon^2}$$

$$\leq 4 \int g_{2\varepsilon} \, d\nu + 4 \frac{\sigma^2}{(2\varepsilon)^2}$$

$$= 4[\nu([-2\varepsilon, 2\varepsilon]^c) + F(2\varepsilon)].$$

The converse inequality follows analogously. □

In general, the following elementary lemma is of great help in the calculations.

LEMMA 5.2. *For $\alpha \in \mathbb{R}$, there exist positive constants $C_1, C_2$ such that, for any $\gamma \geq 0$,*

$$C_1 \frac{e^\gamma - 1 - \gamma - \gamma^2/2}{\gamma} \leq \int_0^1 (e^{\gamma x} - 1 - \gamma x) \frac{dx}{x^\alpha} \leq C_2 \frac{e^\gamma - 1 - \gamma - \gamma^2/2}{\gamma},$$

$$C_1 \frac{e^\gamma - 1 - \gamma}{\gamma} \leq \int_0^1 (e^{\gamma x} - 1) \frac{dx}{x^\alpha} \leq C_2 \frac{e^\gamma - 1 - \gamma}{\gamma},$$

$$C_1 \frac{e^\gamma - 1}{\gamma} \leq \int_0^1 e^{\gamma x} \frac{dx}{x^\alpha} \leq C_2 \frac{e^\gamma - 1}{\gamma},$$

*provided the integral in question converges.*

PROOF. To prove the inequalities, write the exponential as a series, exchange summation and integration (which is possible since everything is absolutely integrable), then integrate term by term. The remaining factor $(n+1)/(n+1-\alpha)$ can be estimated from above and below uniformly in $n$ in the range of the respective sum. □

When $\gamma$ is negative, the last lemma is not valid. Instead, we have the following.

LEMMA 5.3. *Let $\alpha \in \mathbb{R}$. Then, as $\gamma \to -\infty$,*

$$\int_0^1 (e^{\gamma x} - 1 - \gamma x) \frac{dx}{x^\alpha} \approx \begin{cases} \gamma, & \alpha < 2, \\ \gamma \log(-\gamma), & \alpha = 2, \\ -(-\gamma)^{\alpha-1}, & 2 < \alpha < 3, \end{cases}$$

$$\int_0^1 (e^{\gamma x} - 1) \frac{dx}{x^\alpha} \approx \begin{cases} -1, & \alpha < 1, \\ -\log(-\gamma), & \alpha = 1, \\ -(-\gamma)^{\alpha-1}, & 1 < \alpha < 2, \end{cases}$$

$$\int_0^1 e^{\gamma x} \frac{dx}{x^\alpha} \approx (-\gamma)^{\alpha-1} \qquad if \; \alpha < 1.$$

Proof.   Simply substitute $y = -\gamma x$ and calculate the behavior of the integrals at infinity.   □

Furthermore, we need the following fact.

Lemma 5.4.   *Provided the integral in question converges, we have, for $u \in \mathbb{R}$,*

$$\int_0^\varepsilon (e^{ux} - 1 - ux) \frac{dx}{x^\alpha} = \varepsilon^{1-\alpha} \int_0^1 (e^{u\varepsilon x} - 1 - u\varepsilon x) \frac{dx}{x^\alpha},$$

$$\int_0^\varepsilon (e^{ux} - 1) \frac{dx}{x^\alpha} = \varepsilon^{1-\alpha} \int_0^1 (e^{u\varepsilon x} - 1) \frac{dx}{x^\alpha},$$

$$\int_0^\varepsilon e^{ux} \frac{dx}{x^\alpha} = \varepsilon^{1-\alpha} \int_0^1 e^{u\varepsilon x} \frac{dx}{x^\alpha}.$$

We now start with the proofs of the examples. We start with the stable subordinator with drift.

Proof of Example 2.3.   Here, we have to consider

$$\Lambda(u) = u\mu + \int_0^\infty (e^{ux} - 1) \frac{dx}{x^{1+\alpha}}.$$

After subtracting the large jumps, we have

$$\Lambda_\varepsilon(u) = u\mu + \int_0^\varepsilon (e^{ux} - 1) \frac{dx}{x^{1+\alpha}}.$$

We now assume that $u_\varepsilon$ solves

$$0 = \Lambda_\varepsilon'(u_\varepsilon) = \mu + \int_0^\varepsilon e^{u_\varepsilon x} \frac{dx}{x^\alpha} = \mu + \varepsilon^{1-\alpha} \int_0^1 e^{u_\varepsilon \varepsilon x} \frac{dx}{x^\alpha}.$$

Note that $\varepsilon u_\varepsilon$ tends to $\infty$ as $\varepsilon \to 0$. Therefore, setting $u_\varepsilon =: \varepsilon^{-1} \log v$, Lemma 5.2 yields that $v \approx \varepsilon^{\alpha-1} |\log \varepsilon|$. This implies, again using Lemma 5.2, that

$$\nu([-\varepsilon, \varepsilon]^c) \approx \varepsilon^{-\alpha}, \qquad -\Lambda_\varepsilon(u_\varepsilon) \approx \varepsilon^{-1} |\log \varepsilon|,$$

$$\varepsilon |u_\varepsilon| \approx |\log \varepsilon| \quad \text{and} \quad \bar{F}(\varepsilon) \approx \varepsilon^{-1}. \qquad\qquad □$$

In a similar way, we treat the Gamma process with drift.

Proof of Example 2.5.   Here, we have to consider

$$\Lambda(u) = u\mu + \int_0^\infty (e^{ux} - 1) \frac{be^{-x/a}}{x} \, dx.$$

After subtracting the large jumps, we have

$$\Lambda_\varepsilon(u) = u\mu + \int_0^\varepsilon (e^{ux} - 1) \frac{be^{-x/a}}{x} \, dx.$$

Again, we need to consider the solution $u_\varepsilon$ of

$$0 = \Lambda'_\varepsilon(u_\varepsilon) = \mu + \int_0^\varepsilon e^{u_\varepsilon x} b e^{-x/a} \, dx = \mu + \frac{e^{(u_\varepsilon - 1/a)\varepsilon} - 1}{u_\varepsilon - 1/a} b.$$

Noting that $\varepsilon u_\varepsilon \to \infty$ and setting $u_\varepsilon - 1/a =: \varepsilon^{-1} \log v$, we find, using Lemma 5.2, that $v \approx \varepsilon^{-1} |\log \varepsilon|$. This yields, using Lemma 5.2,

$$\nu([-\varepsilon, \varepsilon]^c) \approx |\log \varepsilon|, \qquad -\Lambda_\varepsilon(u_\varepsilon) \approx \varepsilon^{-1} |\log \varepsilon|,$$

$$\varepsilon |u_\varepsilon| \approx |\log \varepsilon| \quad \text{and} \quad \bar{F}(\varepsilon) \approx \varepsilon^{-1}. \qquad \square$$

We now come to the proofs for the polynomial Lévy measure. We distinguish between the asymmetric cases $\alpha_1 < 1$, $\alpha_1 = 1$ and $\alpha_1 > 1$, and the symmetric case.

PROOF OF COROLLARIES 2.7–2.9 AND 2.11. In the symmetric case where $\alpha_1 = \alpha_2$, $C_1 = C_2$ and $b = 0$, the result follows immediately from Corollary 2.1. The remaining cases are treated separately.

*Case* 1: $\alpha_1 \geq \alpha_2$, $\alpha_1 > 1$.

We have

$$
\begin{aligned}
(24) \quad \Lambda_\varepsilon(u) = {}& \int_0^\varepsilon (e^{ux} - 1 - ux) \frac{C_1 \, dx}{x^{1+\alpha_1}} + \int_{-\varepsilon}^0 (e^{ux} - 1 - ux) \frac{C_2 \, dx}{(-x)^{1+\alpha_2}} \\
& + \left( b - \int_\varepsilon^1 x \frac{C_1 \, dx}{x^{1+\alpha_1}} - \int_{-1}^{-\varepsilon} x \frac{C_2 \, dx}{(-x)^{1+\alpha_2}} \right) u.
\end{aligned}
$$

After differentiating with respect to $u$, we obtain

$$
\begin{aligned}
(25) \quad \Lambda'_\varepsilon(u) = {}& \int_0^\varepsilon (e^{ux} - 1) \frac{C_1 \, dx}{x^{\alpha_1}} - \int_0^\varepsilon (e^{-ux} - 1) \frac{C_2 \, dx}{x^{\alpha_2}} \\
& + b - \int_\varepsilon^1 \frac{C_1 \, dx}{x^{\alpha_1}} + \int_\varepsilon^1 \frac{C_2 \, dx}{x^{\alpha_2}}.
\end{aligned}
$$

Note that

$$\Lambda'_\varepsilon(0) = b - \int_\varepsilon^1 \frac{C_1 \, dx}{x^{\alpha_1}} + \int_\varepsilon^1 \frac{C_2 \, dx}{x^{\alpha_2}} \approx -\varepsilon^{1-\alpha_1}$$

if $\alpha_1 > \alpha_2$ or $C_1 > C_2$. Otherwise, one has $\alpha_1 = \alpha_2$ and $C_1 = C_2$, which implies that $\Lambda'_\varepsilon(0) = b$. In the former case, $u_\varepsilon$ is positive for all sufficiently small $\varepsilon > 0$ and, by Lemma 5.4,

$$\varepsilon^{1-\alpha_1} \approx \varepsilon^{1-\alpha_1} \int_0^1 (e^{\varepsilon u_\varepsilon x} - 1) \frac{C_1 \, dx}{x^{\alpha_1}} - \varepsilon^{1-\alpha_2} \int_0^1 (e^{-\varepsilon u_\varepsilon x} - 1) \frac{C_2 \, dx}{x^{\alpha_2}}$$

so that $\varepsilon u_\varepsilon \approx 1$ and, by Lemma 5.2,

$$\nu([-\varepsilon, \varepsilon]^c) \approx \varepsilon^{-\alpha_1}, \qquad -\Lambda(u_\varepsilon) \approx \varepsilon^{-\alpha_1} \quad \text{and} \quad \bar{F}(\varepsilon) \approx \varepsilon^{-\alpha_1}.$$

In the latter case, we obtain

$$\varepsilon^{1-\alpha_1} \int_{-1}^{1} (e^{\varepsilon u_\varepsilon x} - 1) \frac{dx}{x^{\alpha_1}} = \mathcal{O}(1)$$

so that $\varepsilon u_\varepsilon$ tends to zero and

$$\nu([-\varepsilon, \varepsilon]^c) \approx \varepsilon^{-\alpha_1}, \qquad -\Lambda_\varepsilon(u_\varepsilon) = \mathcal{O}(\varepsilon^{-\alpha_1}) \quad \text{and} \quad \bar{F}(\varepsilon) = \mathcal{O}(\varepsilon^{-\alpha_1}).$$

*Case* 2: $\alpha_1 = 1 \geq \alpha_2$ with $\alpha_2 < 1$ or $C_1 > C_2$.

Essentially, we proceed as in the proof of Case 1. Note that (26) is valid, but that

$$\Lambda'_\varepsilon(0) = b - \int_\varepsilon^1 \frac{C_1 \, dx}{x} + \int_\varepsilon^1 \frac{C_2 \, dx}{x^{\alpha_2}} \approx -|\log \varepsilon| \to -\infty.$$

Thus, $u_\varepsilon$ is again positive for sufficiently small $\varepsilon$ and

$$|\log \varepsilon| \approx \int_0^\varepsilon (e^{u_\varepsilon x} - 1) \frac{C_1 \, dx}{x} - \int_0^\varepsilon (e^{u_\varepsilon x} - 1) \frac{C_2 \, dx}{x^{\alpha_2}}$$
$$= \int_0^1 (e^{\varepsilon u_\varepsilon x} - 1) \frac{C_1 \, dx}{x} - \varepsilon^{1-\alpha_2} \int_0^1 (e^{\varepsilon u_\varepsilon x} - 1) \frac{C_2 \, dx}{x^{\alpha_2}}.$$

Now, $\varepsilon u_\varepsilon$ tends to infinity and we use Lemmas 5.2 and 5.3 to deduce that

$$|\log \varepsilon| \approx \frac{e^{\varepsilon u_\varepsilon}}{\varepsilon u_\varepsilon}.$$

Setting $\log v := \varepsilon u_\varepsilon$, we conclude that $v \approx |\log \varepsilon| \log(|\log \varepsilon|)$. This yields, again using Lemmas 5.2 and 5.3, that

$$\nu([-\varepsilon, \varepsilon]^c) \approx \varepsilon^{-1}, \qquad -\Lambda_\varepsilon(u_\varepsilon) \approx \varepsilon^{-1} |\log \varepsilon| \log(|\log \varepsilon|),$$
$$\varepsilon |u_\varepsilon| \approx \log(|\log \varepsilon|) \quad \text{and} \quad \bar{F}(\varepsilon) \approx \varepsilon^{-1} |\log \varepsilon|.$$

*Case* 3: $1 = \alpha_1 = \alpha_2$, $C_1 = C_2$ and $b \neq 0$ (without loss of generality, $b < 0$). Now,

$$\Lambda'_\varepsilon(0) = b$$

and thus $u_\varepsilon$ is positive with

$$\int_{-\varepsilon}^{\varepsilon} (e^{u_\varepsilon x} - 1) \frac{dx}{x^{\alpha_1}} \approx 1,$$

which gives $\varepsilon u_\varepsilon \approx 1$, leading to

$$\nu([-\varepsilon, \varepsilon]^c) \approx \varepsilon^{-1}, \qquad -\Lambda_\varepsilon(u_\varepsilon) \approx \varepsilon^{-1} \quad \text{and} \quad \bar{F}(\varepsilon) \approx \varepsilon^{-1}.$$

*Case* 4: $1 > \alpha_1 \vee \alpha_2$ and $c \neq 0$ (without loss of generality, $c < 0$). Now,

$$\Lambda_\varepsilon'(0) = b - \int_\varepsilon^1 \frac{C_1\, dx}{x^{\alpha_1}} + \int_\varepsilon^1 \frac{C_2\, dx}{x^{\alpha_2}} \to c$$

since $\alpha_1 \vee \alpha_2 < 1$. Therefore,

$$-c \approx \int_0^\varepsilon (e^{ux} - 1) \frac{dx}{x^{\alpha_1}} - \int_0^\varepsilon (e^{-ux} - 1) \frac{dx}{x^{\alpha_2}}.$$

We have $u_\varepsilon > 0$. Set $u_\varepsilon = \varepsilon^{-1} \log v$ and observe that $v \approx \varepsilon^{-(1-\alpha_1)} |\log \varepsilon|$. This gives

$$\nu([-\varepsilon, \varepsilon]^c) \approx \varepsilon^{-\alpha_1 \vee \alpha_2}, \qquad -\Lambda_\varepsilon(u_\varepsilon) \approx \varepsilon^{-1} |\log \varepsilon| \quad \text{and} \quad \bar{F}(\varepsilon) \approx \varepsilon^{-\alpha_1 \vee \alpha_2}.$$

*Case* 5: $1 > \alpha_1 \geq \alpha_2$ and $c = 0$ with $\alpha_1 > \alpha_2$ or $C_1 > C_2$. In this case, $\Lambda_\varepsilon$ simplifies to

$$\Lambda_\varepsilon(u) = \int_0^\varepsilon (e^{ux} - 1) \frac{C_1\, dx}{x^{1+\alpha_1}} + \int_{-\varepsilon}^0 (e^{ux} - 1) \frac{C_2\, dx}{(-x)^{1+\alpha_2}}$$

and we have

$$\Lambda_\varepsilon'(u) = \int_0^\varepsilon e^{ux} \frac{C_1\, dx}{x^{\alpha_1}} - \int_0^\varepsilon e^{-ux} \frac{C_2\, dx}{x^{\alpha_2}}.$$

In particular, $u_\varepsilon$ is negative for sufficiently small $\varepsilon$. Using the above lemmas, one derives that $\varepsilon |u_\varepsilon| \approx 1$ if $\alpha_1 = \alpha_2$ and, otherwise, $\varepsilon |u_\varepsilon| = \log v$ with $v \approx \varepsilon^{\alpha_2 - \alpha_1} |\log \varepsilon|^{\alpha_1}$. In both cases, the term $\nu([-\varepsilon, \varepsilon]^c) \approx \varepsilon^{-\alpha_1}$ is of leading order. □

We finish with the proofs for the remaining examples.

PROOF OF EXAMPLE 2.12. The reasoning is essentially the same as for the Gamma process (Example 2.5) when $\mu \neq 0$. The result follows immediately from Corollary 2.1 for $\mu = 0$. □

PROOF OF EXAMPLE 2.13. The influence of $b_A$ is clear from Theorem 30.1 in [18], so let us assume that $b_A = 0$. We then consider

$$F(\varepsilon) = \varepsilon^{-2} \int_{-\varepsilon}^\varepsilon x^2 \nu(dx) \quad \text{and} \quad N(\varepsilon) := \nu([-\varepsilon, \varepsilon]^c),$$

where $\nu$ is the Lévy measure of $X = B_{A_t}$. It is easy to see (cf. Theorem 30.1 in [18]) that

$$F(\varepsilon) = \varepsilon^{-2} \int_0^\infty \mathbb{E} \mathbb{1}_{\{|B_x| \leq \varepsilon\}} B_x^2 \nu_A(dx)$$

$$(26) \qquad = \varepsilon^{-2} \int_0^\infty \mathbb{E} \mathbb{1}_{\{x\xi^2 \leq \varepsilon^2\}} x \xi^2 \nu_A(dx) = \varepsilon^{-2} \left( \int_0^{\varepsilon^2} \cdots + \int_{\varepsilon^2}^\infty \cdots \right),$$

where $\xi$ is a standard normal random variable. In the first integral in (26), the indicator can be estimated from above by 1 and from below by $\mathbb{1}_{\{\xi^2 \leq 1\}}$. Thus, the first integral in (26) has the order

$$\varepsilon^{-2} \int_0^{\varepsilon^2} x \nu_A(dx) \approx \int_0^{\varepsilon^2} (1 - e^{-x\varepsilon^{-2}}) \nu_A(dx).$$

On the other hand, the second integral in (26) can be estimated from above by

$$\varepsilon^{-2} \int_{\varepsilon^2}^\infty \varepsilon^2 \nu_A(dx).$$

From below, we estimate it by zero. Furthermore,

$$
\begin{aligned}
(27) \quad N(\varepsilon) &= \int_0^\infty \mathbb{P}(|B_x| > \varepsilon) \nu_A(dx) \\
&= \int_0^{\varepsilon^2} \mathbb{P}(x\xi^2 > \varepsilon^2) \nu_A(dx) + \int_{\varepsilon^2}^\infty \mathbb{P}(x\xi^2 > \varepsilon^2) \nu_A(dx).
\end{aligned}
$$

The integrand in the first term in (27) is of order less than

$$e^{-\varepsilon^2/(2x)} \preceq x\varepsilon^{-2} \approx 1 - e^{-x\varepsilon^{-2}}.$$

The integrand in the second term in (27) is bounded from above and below.

Putting all the pieces together, we obtain that

$$
\begin{aligned}
F(\varepsilon) + N(\varepsilon) &\approx \int_0^{\varepsilon^2} (1 - e^{-x\varepsilon^{-2}}) \nu_A(dx) + \int_{\varepsilon^2}^\infty \nu_A(dx) \\
&\approx \int_0^\infty (1 - e^{-x\varepsilon^{-2}}) \nu_A(dx) = \Phi(\varepsilon^{-2}),
\end{aligned}
$$

as required.   □

## APPENDIX

In the examples of Section 2, the term $\varepsilon|u_\varepsilon|$ (appearing in Theorem 1.5) does not affect the small deviation order. We now show that under weak assumptions, one can always neglect the term $\varepsilon|u_\varepsilon|$ since it is of lower order than $-\Lambda_\varepsilon(u_\varepsilon)$. It is still an open question as to whether $\varepsilon|u_\varepsilon|$ is negligible in all cases.

LEMMA A.1.  *Let us assume that we are in the situation of Theorem 1.5.*

(a)  *If $b_\varepsilon \leq 0$, then*

$$(\varepsilon|u_\varepsilon|)^2 \leq -2\left(\varepsilon^{-2} \int_0^\varepsilon x^2 \nu(dx)\right)^{-1} \Lambda_\varepsilon(u_\varepsilon).$$

(b) *If*

$$\int_{-1}^{1} |x| \nu(dx) < \infty \tag{28}$$

*and the effective drift $c \neq 0$, then we have*

$$|cu_\varepsilon| \leq -4\Lambda_\varepsilon(u_\varepsilon)$$

*for all sufficiently small $\varepsilon$.*

REMARK A.2. The assertion implies the negligibility of $\varepsilon|u_\varepsilon|$ in most cases: if $\nu$ has mass infinity, then the term $\varepsilon^{-2} \int_0^\varepsilon x^2 \nu(dx)$ is typically bounded away from zero so that, by assertion (a), $\varepsilon|u_\varepsilon| \leq \mathrm{const} + o(|\Lambda_\varepsilon(u_\varepsilon)|)$; if, on the other hand, $\nu$ has finite mass and $c \neq 0$, then assertion (b) implies that $\varepsilon|u_\varepsilon| = o(|\Lambda_\varepsilon(u_\varepsilon)|)$.

PROOF OF LEMMA A.1. Note that

$$\Lambda_\varepsilon(u) = b_\varepsilon u + \int_{-\varepsilon}^{\varepsilon} (e^{ux} - 1 - ux) \nu(dx),$$

$$b_\varepsilon = b - \int_{-1}^{-\varepsilon} x \nu(dx) - \int_\varepsilon^1 x \nu(dx)$$

and

$$0 = \Lambda_\varepsilon'(u_\varepsilon) = b_\varepsilon + \int_{-\varepsilon}^{\varepsilon} (e^{u_\varepsilon x} - 1) x \nu(dx),$$

which implies that we can express $b_\varepsilon$ in terms of $u_\varepsilon$:

$$b_\varepsilon = - \int_{-\varepsilon}^{\varepsilon} (e^{u_\varepsilon x} - 1) x \nu(dx).$$

This shows that

$$\Lambda_\varepsilon(u_\varepsilon) = - \int_{-\varepsilon}^{\varepsilon} (e^{u_\varepsilon x} - 1) x \nu(dx) u_\varepsilon + \int_{-\varepsilon}^{\varepsilon} (e^{u_\varepsilon x} - 1 - u_\varepsilon x) \nu(dx)$$

$$= \int_{-\varepsilon}^{\varepsilon} (e^{u_\varepsilon x} (1 - u_\varepsilon x) - 1) \nu(dx). \tag{29}$$

*Statement* (a): Since $b_\varepsilon < 0$, we have $u_\varepsilon > 0$. We obtain from (29) and the observation $e^z(1 - z) - 1 \leq 0$ for all real $z$ that

$$\Lambda_\varepsilon(u_\varepsilon) \leq \int_0^\varepsilon (e^{u_\varepsilon x}(1 - u_\varepsilon x) - 1) \nu(dx)$$

$$\leq -\frac{1}{2} \int_0^\varepsilon u_\varepsilon x (e^{u_\varepsilon x} - 1) \nu(dx), \tag{30}$$

where the last step follows from $e^z(1-z) - 1 \leq -\frac{1}{2}z(e^z - 1)$, which holds for all $z \geq 0$. Using $e^z - 1 \geq z$, we can estimate the last term from above by

$$-\frac{1}{2}\int_0^\varepsilon (u_\varepsilon x)^2 \nu(dx) = -\frac{1}{2}(u_\varepsilon \varepsilon)^2 \varepsilon^{-2} \int_0^\varepsilon x^2 \nu(dx),$$

as asserted.

*Statement* (b): Assume that $\varepsilon$ is sufficiently small so that

$$(31) \qquad\qquad \int_{-\varepsilon}^\varepsilon |x| \nu(dx) \leq \frac{|c|}{2}.$$

By (28), we get

$$\Lambda_\varepsilon(u) = cu + \int_{-\varepsilon}^\varepsilon (e^{ux} - 1)\nu(dx)$$

and

$$0 = \Lambda_\varepsilon'(u_\varepsilon) = c + \int_{-\varepsilon}^\varepsilon e^{u_\varepsilon x} x \nu(dx).$$

Without loss of generality, we assume that $c < 0$ and thus $u_\varepsilon > 0$. We estimate

$$(32) \qquad\qquad -c = \int_{-\varepsilon}^\varepsilon e^{u_\varepsilon x} x \nu(dx) \leq \int_0^\varepsilon e^{u_\varepsilon x} x \nu(dx)$$

and conclude, as in (30), that

$$\Lambda_\varepsilon(u_\varepsilon) = cu_\varepsilon + \int_{-\varepsilon}^\varepsilon (e^{u_\varepsilon x} - 1)\nu(dx) \leq -\frac{1}{2}\int_0^\varepsilon u_\varepsilon x(e^{u_\varepsilon x} - 1)\nu(dx),$$

$$= -\frac{u_\varepsilon}{2}\left(\int_0^\varepsilon x e^{u_\varepsilon x} \nu(dx) - \int_0^\varepsilon x \nu(dx)\right) \leq -\frac{u_\varepsilon}{2}\left(-c - \int_0^\varepsilon x \nu(dx)\right).$$

Here, we used (32) in the last step. To estimate the last expression, we use (31), which already yields the assertion. $\square$

**Acknowledgment.** The authors would like to thank Mikhail Lifshits for various valuable comments.

Institut für Mathematik, MA 7-5
Fakultät II
Technische Universität Berlin
Strasse des 17. Juni 136
10623 Berlin
Germany
E-mails: aurzada@math.tu-berlin.de
          dereich@math.tu-berlin.de